\documentclass[10pt,aps,prb,reqno]{amsart}
\usepackage{amsmath}
\usepackage{amssymb}

\newtheorem{theorem}{Theorem}[section]
\newtheorem{remark}{Remark}[section]
\newtheorem{lemma}[theorem]{Lemma}
\newtheorem{proposition}[theorem]{Proposition}

\begin{document}
\title[Two-dimensional MHD system]{Remarks on the global regularity of two-dimensional magnetohydrodynamics system with zero dissipation}
\author{Kazuo Yamazaki}  
\date{}
\thanks{The author expresses gratitude to Professor Jiahong Wu for suggesting this direction of research and Professor David Ullrich for his teaching. The author also expresses gratitude to the referees for helpful suggestions. 
}

\maketitle 

\begin{abstract}
We study the two-dimensional magnetohydrodynamics system with generalized dissipation and diffusion in terms of fractional Laplacians. It is known that the classical magnetohydrodynamics system with full Laplacians in both dissipation and diffusion terms admits a unique global strong solution pair. Making use of the special structure of the system in the two-dimensional case, we show in particular that the solution pair remains smooth when we have zero dissipation but only magnetic diffusion with its power of the fractional Laplacian $\beta > \frac{3}{2}$. 

\vspace{5mm}

\end{abstract}
\footnote{2000MSC : 35B65, 35Q35, 35Q86}
\footnote{Department of Mathematics, Oklahoma State University, 401 Mathematical Sciences, Stillwater, OK 74078, USA}

\section{Introduction and statement of results}

We study the following magnetohydrodynamics (MHD) system: 

\begin{equation}
\begin{cases}
\frac{\partial u}{\partial t} + (u\cdot\nabla) u - (b\cdot\nabla)b + \nabla \pi + \nu \Lambda^{2\alpha}u = 0 \\
\frac{\partial b}{\partial t} + (u\cdot\nabla) b - (b\cdot\nabla)u + \eta \Lambda^{2\beta} b = 0\\
\nabla\cdot u = \nabla\cdot b = 0, \hspace{3mm} (u,b)(x,0) = (u_{0}, b_{0})(x)
\end{cases}
\end{equation}

where $u: \mathbb{R}^{N}\times \mathbb{R}^{+}\mapsto \mathbb{R}^{N}$ represents the velocity vector field, $b: \mathbb{R}^{N}\times \mathbb{R}^{+}\mapsto \mathbb{R}^{N}$ the magnetic vector field, $\pi: \mathbb{R}^{N}\times \mathbb{R}^{+}\mapsto \mathbb{R}$ the pressure scalar field and $\nu, \eta \geq 0$ are the kinematic viscosity and diffusivity constants respectively. We also let $\hat{f}(\xi)$ denote the Fourier transform of $f$; i.e.

\begin{equation*}
 \hat{f}(\xi) = \int_{\mathbb{R}^{N}}f(x) e^{-ix\cdot\xi}dx
 \end{equation*} 

and defined a fractional Laplacian operator $\Lambda^{2\gamma}$ with $\gamma \in \mathbb{R}$ to have the Fourier symbol of $\lvert\xi\rvert^{2\gamma}$; that is, 

\begin{equation*}
\widehat{\Lambda^{2\gamma}f}(\xi) = \lvert \xi\rvert^{2\gamma} \hat{f}(\xi)
\end{equation*}

In case $N = 2, 3, \nu, \eta > 0, \alpha = \beta = 1$, the MHD system possesses at least one global $L^{2}$ weak solution for any initial data pair $(u_{0}, b_{0}) \in L^{2}(\mathbb{R}^{N})\times L^{2}(\mathbb{R}^{N})$; in case $N =2$, in fact the solution is unique (cf. [20]). 

In order to discuss the previous results on strong solutions and better understand the importance of the lower bounds for the two parameters $\alpha, \beta > 0$ when $\nu, \eta > 0$, let us recall the notion of criticality in a simple setting. Firstly, it can be shown that the solution pair to (1) with $\alpha = \beta = \gamma$ has the rescaling properties that if $(u(x,t), b(x,t))$ solves the system, then so does $\left(u_{\lambda}(x,t), b_{\lambda}(x,t)\right)$ with $\lambda \in \mathbb{R}^{+}$ where

\begin{equation*}
u_{\lambda}(x,t) = \lambda^{2\gamma -1} u(\lambda x, \lambda^{2\gamma} t), \hspace{3mm} b_{\lambda}(x,t) = \lambda^{2\gamma -1} b(\lambda x, \lambda^{2\gamma}t), \hspace{3mm} \gamma \in \mathbb{R}^{+}
\end{equation*}

As we show in (4), the solution pair $(u,b)$ to (1) has the global bounds on the $L^{2}$-norm and it can be shown that $\gamma = \frac{1}{2} + \frac{N}{4}$ implies 

\begin{equation*}
\lVert u_{\lambda}(\cdot, t)\rVert_{L^{2}(\mathbb{R}^{N})} = \lVert u(\cdot, \lambda^{2\gamma} t)\rVert_{L^{2}(\mathbb{R}^{N})}, \hspace{3mm} \lVert b_{\lambda}(\cdot, t)\rVert_{L^{2}(\mathbb{R}^{N})} = \lVert b(\cdot, \lambda^{2\gamma} t)\rVert_{L^{2}(\mathbb{R}^{N})}
\end{equation*}

With this in mind, we call the case $\nu, \eta > 0, \alpha \geq \frac{1}{2} + \frac{N}{4}, \beta \geq \frac{1}{2} + \frac{N}{4}$ the critical case and in such a case, the existence of the unique global strong solution pair has been shown (cf. [25]). 

Some numerical analysis results (e.g. [11], [19]) indicate more dominant role played by the velocity vector field in preserving the regularity of the solution pair. Moreover, starting from the works of [12] and [31], we have also seen various regularity criteria of the MHD system in terms of only the velocity vector field (e.g. [2], [6], [8], [10], [13], [24], [27], [30]). This is largely due to the fact that upon taking $H^{1}$-estimates of $u$ and $b$, every non-linear term involves $u$ while not necessarily $b$. With this in mind, following the work of [21], the author in [23] showed that even in logarithmically super-critical case the system (1) still admits a unique global strong solution pair. That is, the author replaced the dissipative term of $\nu\Lambda^{2\alpha} u$ and the diffusive term of $\eta \Lambda^{2\beta} b$ by $\nu\mathcal{L}_{1}^{2}u$ and $\eta \mathcal{L}_{2}^{2}b$ respectively where $\mathcal{L}_{i}, i = 1, 2$ are defined to have the Fourier symbols of $m_{i}(\xi), i = 1, 2$ satisfying the following lower bounds: 

\begin{equation*}
\widehat{\mathcal{L}_{1}u}(\xi) = m_{1}(\xi) \hat{u}(\xi), \hspace{3mm} \widehat{\mathcal{L}_{2}b}(\xi) = m_{2}(\xi) \hat{b}(\xi)
\end{equation*}

and 

\begin{equation*}
m_{1}(\xi) \geq \frac{\lvert \xi\rvert^{\alpha}}{g_{1}(\xi)}, \hspace{3mm} m_{2}(\xi) \geq \frac{\lvert\xi\rvert^{\beta}}{g_{2}(\xi)}, \hspace{3mm} \alpha \geq \frac{1}{2} + \frac{N}{4}, \beta > 0, \alpha + \beta \geq 1 + \frac{N}{2}
\end{equation*}

with $g_{i} \geq 1, i = 1 ,2$ being radially symmetric, non-decreasing functions. 

The endpoint case $\nu > 0, \eta = 0, \alpha = 1 + \frac{N}{2}$ was also completed recently in [28] (cf. also [26] for further generalization). 

On the other hand, in case $N =2$, it is well-known that the Euler equation, the Navier-Stokes system with no dissipation, admits a unique global strong solution. This is due to the fact that upon taking a curl, the vorticity becomes a conserved quantity. In the case of the MHD system, upon taking a curl and $L^{2}$-estimate of the resulting system, every non-linear term has $b$ involved. Exploiting this observation and divergence-free conditions, the authors in [3] showed that in case $N = 2$, full Laplacians in both dissipation and magnetic diffusion are not necessary for the solution to remain smooth; rather, only a mix of partial dissipation and diffusion in the order of two derivatives suffices. In this paper we make further observation in case $N = 2$:

\begin{theorem}
Let $N = 2, \nu = 0, \eta > 0, \alpha = 0, \beta > \frac{3}{2}$. Then for all initial data pair $(u_{0}, b_{0}) \in H^{s}(\mathbb{R}^{2}) \times H^{s}(\mathbb{R}^{2}), s \geq 1 + 2\beta$, there exists a unique global strong solution pair $(u,b)$ to (1) such that 

\begin{eqnarray*}
&&u \in C([0,\infty); H^{s}(\mathbb{R}^{2}))\\
&&b \in C([0,\infty); H^{s}(\mathbb{R}^{2}))\cap L^{2}([0,\infty); H^{s+\beta}(\mathbb{R}^{2}))
\end{eqnarray*}

\end{theorem}

\begin{theorem}
Let $N = 2, \nu, \eta > 0, \alpha \in \left(0, \frac{1}{2}\right), \beta \in \left(\frac{5}{4}, \frac{3}{2}\right]$ such that $\alpha + 2\beta > 3$. Then for all initial data pair $(u_{0}, b_{0}) \in H^{s}(\mathbb{R}^{2}) \times H^{s}(\mathbb{R}^{2}), s \geq 1+2\beta$, there exists a unique global strong solution pair $(u,b)$ to (1) such that 

\begin{eqnarray*}
&&u \in C([0,\infty); H^{s}(\mathbb{R}^{2}))\cap L^{2}([0,\infty); H^{s+\alpha} (\mathbb{R}^{2}))\\
&&b \in C([0,\infty); H^{s}(\mathbb{R}^{2}))\cap L^{2}([0,\infty); H^{s+\beta}(\mathbb{R}^{2}))
\end{eqnarray*}
\end{theorem}

\begin{remark}
\begin{enumerate}
\item Our proof was inspired partially from the work of [3], [4] and [6]. We note that making use of the structure of the partial differential equation has proven to be useful in other cases as well (e.g. [29]).
\item While this paper was being prepared, the work by [22] appeared. In their work, it is shown that in particular if $\alpha = 0$, then $\beta > 2$ is required (See Theorem 1 and Remark 1 of [22]) while our Theorem 1.1 shows that $\beta > \frac{3}{2}$ suffices. We also independently obtained Theorem 5.1; this is no longer a new result and thus we placed this in the Appendix because its proof is immediate and very simple. The hypothesis of Theorem 5.1 allows $\alpha \geq \frac{1}{2}$ rather than $\alpha = 0$ as in Theorem 1.1. As will be discussed, a complete lack of dissipation makes the analysis significantly more difficult in the latter case.   
\item There are ways to obtain different initial regularity in various space of functions; we chose to state the above for simplicity. We also refer readers to [3] where the authors considered the case $N = 2, \nu = 0, \eta > 0, \beta = 1$ and showed the existence of weak solution pair and regularity criteria for its global regularity and uniqueness (cf. also [25]). 
\item To extend such a type of result to higher dimension, it seems to require a new idea. As indicated in the work of [23] and [28], in higher dimension, dissipation seems to be crucial in preserving the regularity of the solution pair. 
\end{enumerate}

\end{remark}

In the Preliminary section, let us briefly set up notations and state key lemmas; thereafter, we prove our theorems. 

\section{Preliminary}

Let us denote a constant that depends on $a, b$ by $c(a,b)$ and also  denote curl $u$ by $w$ and similarly $\text{curl } b = j$. We also  denote partial derivatives as follows:

\begin{equation*}
\frac{\partial}{\partial t} = \partial_{t}, \hspace{5mm} \frac{\partial}{\partial x} = \partial_{1}, \hspace{5mm} \frac{\partial}{\partial y} = \partial_{2}
\end{equation*}
 
For simplicity we also set 

\begin{eqnarray}
&&X(t) = \lVert w(t)\rVert_{L^{2}}^{2} + \lVert j(t)\rVert_{L^{2}}^{2}
\end{eqnarray}
 
We use the following well-known inequality (cf. [5]): 

\begin{lemma}
(cf. [5]) Let $f$ be divergence-free vector field such that $\nabla f \in L^{p}, p \in (1,\infty)$. Then there exists a constant $c > 0$ such that

\begin{equation*}
\lVert \nabla f\rVert_{L^{p}} \leq c\frac{p^{2}}{p-1}\lVert \text{curl }  f\rVert_{L^{p}}
\end{equation*}
\end{lemma}

We will use the following commutator estimate:

\begin{lemma}
(cf. [15]) Let $f,g$ be smooth such that $\nabla f \in L^{p_{1}}, \Lambda^{s-1}g \in L^{p_{2}}, \Lambda^{s}f \in L^{p_{3}}, g \in L^{p_{4}}, p \in (1,\infty), \frac{1}{p} = \frac{1}{p_{1}}+\frac{1}{p_{2}} = \frac{1}{p_{3}} + \frac{1}{p_{4}}, p_{2}, p_{3} \in (1,\infty), s > 0.$ Then there exists a constant $c > 0$ such that

\begin{equation*}
\lVert \Lambda^{s}(fg) - f\Lambda^{s}g\rVert_{L^{p}} \leq c(\lVert \nabla f\rVert_{L^{p_{1}}}\lVert \Lambda^{s-1}g\rVert_{L^{p_{2}}} + \lVert \Lambda^{s}f\rVert_{L^{p_{3}}}\lVert g\rVert_{L^{p_{4}}})
\end{equation*}

\end{lemma}

The following logarithmic inequality starting from the works of [1] and [9] has been proven to be useful:

\begin{lemma}
Let $f \in L^{2}(\mathbb{R}^{2})\cap H^{s}(\mathbb{R}^{2}), s > 2$ and f be a divergence-free vector field that satisfies $\text{ curl } f \in L^{\infty}(\mathbb{R}^{2})$. Then there exists a constant $c > 0$ such that 

\begin{equation*}
 \lVert \nabla f\rVert_{L^{\infty}} \leq c\left(\lVert f\rVert_{L^{2}} + \lVert \text{curl } f\rVert_{L^{\infty}} \log_{2}(2+ \lVert f\rVert_{H^{s}}) + 1\right)
 \end{equation*} 

\end{lemma}

For the paper to be self-contained, we sketch its proof in the Appendix. 

\begin{lemma}
(cf. [7], [14])
For any $\alpha \in [0, 1], x \in \mathbb{R}^{N}, \mathbb{T}^{N}, N \in \mathbb{N}$ and $f, \Lambda^{2\alpha}f \in L^{p}, p \geq 2$, 

\begin{equation*}
2\int \lvert\Lambda^{\alpha}(f^{\frac{p}{2}})\rvert^{2} \leq p \int\lvert f\rvert^{p-2}f\Lambda^{2\alpha}f
\end{equation*}

\end{lemma}

Finally, the following product estimate appeared in [16], [17] and [18]: 

\begin{lemma}
Let $\sigma_{1}, \sigma_{2} < 1, \sigma_{1} + \sigma_{2} > 0$. Then there exists a constant $c(\sigma_{1}, \sigma_{2}) > 0$ such that 

\begin{equation*}
\lVert fg\rVert_{\dot{H}^{\sigma_{1} + \sigma_{2} - 1}} \leq c(\sigma_{1}, \sigma_{2}) \lVert f\rVert_{\dot{H}^{\sigma_{1}}}\lVert g\rVert_{\dot{H}^{\sigma_{2}}}
\end{equation*}

for $f \in \dot{H}^{\sigma_{1}}(\mathbb{R}^{2}), g\in \dot{H}^{\sigma_{2}}(\mathbb{R}^{2})$. 

\end{lemma}
We remark that this result may be generalized to $N$-dimension. Let us prove this in the Appendix as well. 

\section{Proof of Theorem 1.1} 

We now work on 

\begin{equation}
\begin{cases}
\partial_{t}u + (u\cdot\nabla)u - (b\cdot\nabla) b + \nabla \pi = 0\\
\partial_{t}b + (u\cdot\nabla) b - (b\cdot\nabla) u + \eta \Lambda^{2\beta}b = 0
\end{cases}
\end{equation}

and assume $\beta \in \left(\frac{3}{2}, 2\right)$ as the case $\beta \geq 2$ may be done after a slight modification. 

Firstly, taking $L^{2}$-inner products of the first equation with $u$ and the second equation with $b$, we obtain in sum 

\begin{eqnarray*}
&&\frac{1}{2}\partial_{t}(\lVert u \rVert_{L^{2}}^{2} + \lVert b\rVert_{L^{2}}^{2}) + \eta \lVert \Lambda ^{\beta}b\rVert_{L^{2}}^{2}\\
&=& -\int (u\cdot\nabla)u\cdot u + (u\cdot\nabla)b\cdot b + \int b\cdot (\nabla b \cdot u + \nabla u \cdot b) - \int \nabla \pi \cdot u
\end{eqnarray*}

and hence using incompressibility conditions and integrating in time, we obtain 

\begin{equation}
\sup_{t\in [0,T]} \left(\lVert u(t)\rVert_{L^{2}}^{2} + \lVert b(t)\rVert_{L^{2}}^{2}\right) + 2\eta\int_{0}^{T} \lVert \Lambda^{\beta} b\rVert_{L^{2}}^{2} d\tau \leq c(u_{0}, b_{0}, T)
\end{equation}

Using this, we obtain the following proposition: 

\begin{proposition}
Let $N = 2, \nu = 0, \eta > 0, \alpha =0, \beta > \frac{3}{2}$. Then for any  solution pair $(u,b)$ to (1) in $[0,T]$, there exists a constant $c(u_{0}, b_{0}, T) > 0$ such that   

\begin{equation*}
\sup_{t\in [0,T]}\left(\lVert w(t)\rVert_{L^{2}}^{2} + \lVert j(t)\rVert_{L^{2}}^{2}\right) + \eta\int_{0}^{T}\lVert \Lambda^{\beta} j\rVert_{L^{2}}^{2} d\tau \leq c(u_{0}, b_{0}, T)
\end{equation*}

\end{proposition}

\proof{

We take a curl on the system (3) to obtain

\begin{eqnarray}
&&\partial_{t} w = - (u\cdot\nabla)w + (b\cdot\nabla)j\\
&&\partial_{t} j + \eta \Lambda^{2\beta}j = - (u\cdot\nabla)j  +  (b\cdot\nabla)w + 2[\partial_{1}b_{1}(\partial_{1}u_{2} + \partial_{2}u_{1}) - \partial_{1}u_{1}(\partial_{1}b_{2} + \partial_{2}b_{1})]\nonumber
\end{eqnarray}

We take $L^{2}$-inner products of (5) with $w$ and $j$ respectively and sum to obtain due to the incompressibility conditions  

\begin{equation}
\frac{1}{2} \partial_{t}X(t) + \eta \lVert \Lambda^{\beta} j\rVert_{L^{2}}^{2} = 2\int [\partial_{1}b_{1}(\partial_{1}u_{2} + \partial_{2}u_{1}) - \partial_{1}u_{1}(\partial_{1}b_{2} + \partial_{2}b_{1})]j
\end{equation}

Now we estimate from (6) as follows: 

\begin{equation*}
\frac{1}{2} \partial_{t}X(t) + \eta \lVert \Lambda^{\beta} j\rVert_{L^{2}}^{2} \leq c\int \lvert \nabla b\rvert \lvert \nabla u\rvert \lvert j\rvert \leq c\lVert \nabla b\rVert_{L^{\frac{2}{\beta - 1}}} \lVert \nabla u\rVert_{L^{2}} \lVert j\rVert_{L^{\frac{2}{2-\beta}}}
\end{equation*}

by H$\ddot{o}$lder's inequalities. A Gagliardo-Nirenberg inequality of 

\begin{equation*}
 \lVert \nabla b\rVert_{L^{\frac{2}{\beta - 1}}} \leq c\lVert b\rVert_{L^{2}}^{\frac{2(\beta - 1)}{1+\beta}}\lVert \Lambda^{1+\beta}b\rVert_{L^{2}}^{1-\frac{2(\beta - 1)}{1+\beta}}
 \end{equation*} 

Lemma 2.1 and Sobolev embedding of $\dot{H}^{\beta - 1}(\mathbb{R}^{2})\hookrightarrow L^{\frac{2}{2-\beta}}(\mathbb{R}^{2})$ lead to the bound of right hand side by 

\begin{equation*}
c\lVert \nabla b\rVert_{L^{\frac{2}{\beta - 1}}} \lVert \nabla u\rVert_{L^{2}} \lVert j\rVert_{L^{\frac{2}{2-\beta}}} \leq c\lVert b\rVert_{L^{2}}^{\frac{2(\beta - 1)}{1+\beta}}\lVert \Lambda^{1+\beta}b\rVert_{L^{2}}^{1-\frac{2(\beta - 1)}{1+\beta}}\lVert w\rVert_{L^{2}} \lVert \Lambda^{\beta -1}j\rVert_{L^{2}}
\end{equation*}

Using 

\begin{equation*}
\sup_{t \in [0,T]} \lVert b(t)\rVert_{L^{2}} \leq c(u_{0}, b_{0}, T)
\end{equation*}

from (4), Lemma 2.1 and Young's inequalities, we further bound by

\begin{eqnarray*}
&&c\lVert \Lambda^{\beta} j\rVert_{L^{2}}^{\frac{3-\beta}{1+\beta}}\lVert w\rVert_{L^{2}} \lVert \Lambda^{\beta} b\rVert_{L^{2}}\\
&\leq& \frac{\eta}{2} \lVert \Lambda^{\beta} j\rVert_{L^{2}}^{2} + c\lVert w\rVert_{L^{2}}^{\frac{2(1+\beta)}{3\beta - 1}}\lVert \Lambda^{\beta}b\rVert_{L^{2}}^{\frac{2(1+\beta)}{3\beta - 1}}\\
&\leq& \frac{\eta}{2} \lVert \Lambda^{\beta} j\rVert_{L^{2}}^{2} + + c(1+ X(t))(1+ \lVert \Lambda^{\beta} b\rVert_{L^{2}}^{2})
\end{eqnarray*}

Using this bound, absorbing the diffusive term, we obtain 

\begin{equation*}
\partial_{t}X(t) \leq c(1+X(t)) (1+ \lVert \Lambda^{\beta} b\rVert_{L^{2}}^{2})
\end{equation*}

which implies 

\begin{equation*}
\sup_{t \in [0,T]} X(t) \leq X(0)e^{c\int_{0}^{T} 1 + \lVert \Lambda^{\beta} b\rVert_{L^{2}}^{2} d\tau } \leq c(u_{0}, b_{0}, T)
\end{equation*}

by (4). It follows that 

\begin{equation*}
\int_{0}^{T} \lVert \Lambda^{\beta} j\rVert_{L^{2}}^{2} d\tau \leq c(u_{0}, b_{0}, T)
\end{equation*}

The proof of Proposition 3.1 is complete. 

}

Next, we first obtain a higher regularity in the magnetic field only making use of diffusivity as much as possible: 

\begin{proposition}
Let $N = 2, \nu = 0, \eta > 0, \alpha = 0, \beta > \frac{3}{2}$. Then for any  solution pair $(u,b)$ to (1) in $[0,T]$, there exists a constant $c(u_{0}, b_{0}, T) > 0$ such that   

\begin{equation*}
\sup_{t\in [0,T]}\lVert \Lambda^{\beta} b(t)\rVert_{L^{2}}^{2} + \eta\int_{0}^{T} \lVert \Lambda^{2\beta} b\rVert_{L^{2}}^{2} d\tau \leq c(u_{0}, b_{0}, T)
\end{equation*}

\end{proposition}

\proof{
Applying $\Lambda^{\beta}$ on the second equation of (3) and taking $L^{2}$-inner products with $\Lambda^{\beta}b$, we have

\begin{equation*}
 \frac{1}{2}\partial_{t} \lVert \Lambda^{\beta} b\rVert_{L^{2}}^{2} + \eta \lVert \Lambda^{2\beta} b\rVert_{L^{2}}^{2} = -\int (u\cdot\nabla)b\cdot\Lambda^{2\beta} b + \int (b\cdot\nabla) u\cdot\Lambda^{2\beta} b
 \end{equation*} 

We bound the right hand side by 

\begin{eqnarray*}
&&\lVert u\rVert_{L^{2}} \lVert \nabla b\rVert_{L^{\infty}} \lVert \Lambda^{2\beta} b\rVert_{L^{2}} + \lVert b\rVert_{L^{\infty}} \lVert \nabla u\rVert_{L^{2}} \lVert \Lambda^{2\beta} b\rVert_{L^{2}}\\
&\leq& c (\lVert \Lambda b\rVert_{L^{2}}^{\frac{\beta - 1}{\beta}}\lVert \Lambda^{1+\beta}b\rVert_{L^{2}}^{\frac{1}{\beta}}\lVert \Lambda^{2\beta} b\rVert_{L^{2}} + \lVert b\rVert_{L^{2}}^{\frac{\beta -1}{\beta}}\lVert\Lambda^{\beta} b\rVert_{L^{2}}^{\frac{1}{\beta}} \lVert w\rVert_{L^{2}} \lVert \Lambda^{2\beta} b\rVert_{L^{2}})
\end{eqnarray*}

due to H$\ddot{o}$lder's inequalities, 

\begin{equation*}
\sup_{t\in [0,T]} \lVert u(t)\rVert_{L^{2}} \leq c(u_{0}, b_{0}, T)
\end{equation*}

from (4), Gagliardo-Nirenberg inequalities of 

\begin{equation}
\lVert f\rVert_{L^{\infty}} \leq c_{0}\lVert f\rVert_{L^{2}}^{\frac{\beta - 1}{\beta}}\lVert \Lambda^{\beta} f\rVert_{L^{2}}^{\frac{1}{\beta}}
\end{equation}

for some constant $c_{0} > 0$ independent of $f$ and Lemma 2.1. Next, the bound on $\lVert b\rVert_{L^{2}}$ from (4), the bound on $\lVert w\rVert_{L^{2}}$ and $\lVert j\rVert_{L^{2}}$ from Proposition 3.1 and Young's inequalities lead to a further bound of 

\begin{eqnarray*}
c(\lVert \Lambda^{1+\beta}b\rVert_{L^{2}}^{\frac{1}{\beta}}\lVert \Lambda^{2\beta} b\rVert_{L^{2}} + \lVert\Lambda^{\beta} b\rVert_{L^{2}}^{\frac{1}{\beta}} \lVert \Lambda^{2\beta} b\rVert_{L^{2}})\leq \frac{\eta}{2} \lVert \Lambda^{2\beta} b\rVert_{L^{2}}^{2} + c(\lVert \Lambda^{1+\beta}b\rVert_{L^{2}}^{\frac{2}{\beta}} + \lVert \Lambda^{\beta}b\rVert_{L^{2}}^{\frac{2}{\beta}})
\end{eqnarray*}

We use this estimate, absorb the diffusive term, rely on Young's inequality and integrate in time to obtain

\begin{equation*}
\sup_{t\in [0,T]} \lVert \Lambda^{\beta} b(t)\rVert_{L^{2}}^{2} + \eta\int_{0}^{T} \lVert \Lambda^{2\beta} b\rVert_{L^{2}}^{2} d\tau \leq \lVert \Lambda^{\beta} b_{0}\rVert_{L^{2}}^{2} + c\int_{0}^{T} \lVert \Lambda^{\beta}j\rVert_{L^{2}}^{2} + \lVert \Lambda^{\beta} b\rVert_{L^{2}}^{2} + 1 d\tau 
\end{equation*}

According to (4) and Proposition 3.1, the time integral on the right hand side is bounded by some constant $c(u_{0}, b_{0}, T)$. This completes the proof of Proposition 3.2. 

}

The higher regularity of the magnetic field from Proposition 3.2 allows us to prove the following proposition: 

\begin{proposition}
Let $N = 2, \nu = 0, \eta > 0, \alpha = 0, \beta > \frac{3}{2}$. Then for any solution pair $(u, b)$ to (1) in $[0,T]$, there exists a constant $c(u_{0}, b_{0}, T) > 0$ such that  

\begin{equation*}
\sup_{t\in [0,T]}\lVert w(t)\rVert_{L^{\infty}}\leq c(u_{0}, b_{0}, T)
\end{equation*}

\end{proposition}

\proof{

We fix $p > 2$, multiply the first equation of (5) with $\lvert w\rvert^{p-2}w$ and integrate in space to estimate by 

\begin{eqnarray*}
\frac{1}{p} \partial_{t} \lVert w\rVert_{L^{p}}^{p} = -\int (u\cdot\nabla) w \lvert w\rvert^{p-2}w + \int (b\cdot\nabla)j \lvert w\rvert^{p-2}w
\leq \lVert b\rVert_{L^{p}} \lVert \nabla j\rVert_{L^{\infty}} \lVert w\rVert_{L^{p}}^{p-1}
\end{eqnarray*}

where we used the divergence-free property of $u$ and H$\ddot{o}$lder's inequality. Dividing by $\lVert w\rVert_{L^{p}}^{p-1}$, we further estimate by 

\begin{eqnarray*}
\partial_{t} \lVert w\rVert_{L^{p}} &\leq& \lVert b\rVert_{L^{p}} \lVert \nabla j\rVert_{L^{\infty}}\\
&\leq& c(p)c \lVert b\rVert_{L^{2}}^{\frac{2+ p(\beta -1)}{\beta p}}\lVert \Lambda^{\beta} b\rVert_{L^{2}}^{1-\frac{2 + p(\beta -1)}{\beta p}}\lVert \Lambda^{\beta} b\rVert_{L^{2}}^{\frac{2\beta -3}{\beta}}\lVert \Lambda^{2\beta} b\rVert_{L^{2}}^{\frac{3-\beta}{\beta}}
\end{eqnarray*}

due to a Gagliardo-Nirenberg inequalities of 

\begin{equation*}
\lVert b\rVert_{L^{p}} \leq c(p) \lVert b\rVert_{L^{2}}^{\frac{2+p(\beta -1)}{\beta p}}\lVert \Lambda^{\beta} b\rVert_{L^{2}}^{1-\frac{2+ p(\beta -1)}{\beta p}}
\end{equation*}

and another which requires $\beta > \frac{3}{2}$: 

\begin{equation*}
\lVert \nabla j\rVert_{L^{\infty}} \leq c \lVert \Lambda^{\beta -1} j\rVert_{L^{2}}^{\frac{2\beta - 3}{\beta}}\lVert \Lambda^{2\beta - 1} j\rVert_{L^{2}}^{\frac{3-\beta}{\beta}}
\end{equation*}

By the bound on $\lVert b\rVert_{L^{2}}$ from (4) and the bound on $\lVert \Lambda^{\beta} b\rVert_{L^{2}}$ from Proposition 3.2, we obtain

\begin{eqnarray*}
\partial_{t} \lVert w\rVert_{L^{p}} 
\leq c(p)c \lVert \Lambda^{2\beta} b\rVert_{L^{2}}^{\frac{3-\beta}{\beta}}
\end{eqnarray*}

Integrating in time over $[0,t]$, by Young's inequality we have 

\begin{equation}
\lVert w(t)\rVert_{L^{p}} \leq \lVert w(0)\rVert_{L^{p}}+ c(p)c\int_{0}^{T} \lVert \Lambda^{2\beta }b\rVert_{L^{2}}^{2} + 1d\tau
\end{equation}

We take limit $p \to \infty$ on (8) and due to the Gagliardo-Nirenberg inequality (7), we obtain

\begin{equation*}
\lVert w(t)\rVert_{L^{\infty}} \leq \lVert w(0)\rVert_{L^{\infty}} + c_{0}c\int_{0}^{T} \lVert \Lambda^{2\beta} b\rVert_{L^{2}}^{2} + 1d\tau 
\end{equation*}

By Proposition 3.2, the right hand side is bounded by $c(u_{0}, b_{0}, T)$. This completes the proof of Proposition 3.3. 

\textit{Proof of Theorem 1.1}

It is well-known that Proposition 3.3 leads to the global regularity of the solution pair $(u,b)$ to (3). We sketch for completeness. 

By Lemma 2.3, the bound on $\lVert u\rVert_{L^{2}}$ from (4) and Proposition 3.3, for any $\gamma > 2$ we have

\begin{eqnarray}
\lVert \nabla u\rVert_{L^{\infty}} &\leq& c\left(\lVert u\rVert_{L^{2}} + \lVert w\rVert_{L^{\infty}} \log_{2}(2+ \lVert u\rVert_{H^{\gamma}}) + 1\right)\\
&\leq& c(u_{0}, b_{0}, T)(\log_{2} (2+ \lVert u\rVert_{H^{\gamma}}) + 1)\nonumber
\end{eqnarray}

Applying $\Lambda^{\gamma}$ on (3), taking $L^{2}$-inner products with $\Lambda^{\gamma}u$ and $\Lambda^{\gamma}b$ respectively we obtain 

\begin{eqnarray}
&&\frac{1}{2}\partial_{t} (\lVert \Lambda^{\gamma}u\rVert_{L^{2}}^{2} + \lVert \Lambda^{\gamma} b\rVert_{L^{2}}^{2}) + \eta \lVert \Lambda^{\gamma + \beta} b\rVert_{L^{2}}^{2}\\
&=& -\int \Lambda^{\gamma}[(u\cdot\nabla)u]\cdot \Lambda^{\gamma}u - u\cdot\nabla \Lambda^{\gamma} u\cdot\Lambda^{\gamma}u - \int \Lambda^{\gamma}[(u\cdot\nabla)b]\cdot \Lambda^{\gamma}b - u\cdot\nabla\Lambda^{\gamma}b\cdot\Lambda^{\gamma}b\nonumber\\
&&+ \int \Lambda^{\gamma}[(b\cdot\nabla)b]\cdot\Lambda^{\gamma}u - b\cdot\nabla \Lambda^{\gamma}b\cdot\Lambda^{\gamma}u +\int \Lambda^{\gamma}[(b\cdot\nabla)u]\cdot\Lambda^{\gamma} b - b\cdot\nabla \Lambda^{\gamma}u\cdot\Lambda^{\gamma}b\nonumber
\end{eqnarray}

because by incompressibility 

\begin{equation*}
\int u\cdot\nabla \Lambda^{\gamma} u\cdot\Lambda^{\gamma}u = \int u\cdot\nabla\Lambda^{\gamma}b\cdot\Lambda^{\gamma}b = 0, \hspace{5mm} 
\int b\cdot\nabla \Lambda^{\gamma}b\cdot\Lambda^{\gamma}u + b\cdot\nabla \Lambda^{\gamma}u\cdot\Lambda^{\gamma}b  = 0
\end{equation*}

By Lemma 2.2 and (9) we obtain 

\begin{eqnarray*}
&&\partial_{t} (\lVert \Lambda^{\gamma}u\rVert_{L^{2}}^{2} + \lVert \Lambda^{\gamma} b\rVert_{L^{2}}^{2}) + 2\eta \lVert \Lambda^{\gamma + \beta} b\rVert_{L^{2}}^{2}\\
&\leq& c (\lVert \nabla u\rVert_{L^{\infty}} + \lVert \nabla b\rVert_{L^{\infty}})(\lVert \Lambda^{\gamma}u\rVert_{L^{2}}^{2} + \lVert \Lambda^{\gamma} b\rVert_{L^{2}}^{2})\\
&\leq& c(u_{0}, b_{0}, T)(\log_{2} (2+ \lVert u\rVert_{H^{\gamma}}  )+1 + \lVert \nabla b\rVert_{L^{2}}^{\frac{2\beta - 2}{2\beta - 1}}\lVert \Lambda^{2\beta}b\rVert_{L^{2}}^{\frac{1}{2\beta -1}})(\lVert \Lambda^{\gamma}u\rVert_{L^{2}}^{2} + \lVert \Lambda^{\gamma} b\rVert_{L^{2}}^{2})\\
&\leq& c(u_{0}, b_{0}, T)(\log_{2} (2+ \lVert u\rVert_{H^{\gamma}}  )+1 + \lVert \Lambda^{2\beta}b\rVert_{L^{2}}^{2})(\lVert \Lambda^{\gamma}u\rVert_{L^{2}}^{2} + \lVert \Lambda^{\gamma} b\rVert_{L^{2}}^{2})\\
&\leq& c(u_{0}, b_{0}, T)(1+ \lVert \Lambda^{2\beta}b\rVert_{L^{2}}^{2})(\log_{2}(2+ \lVert u\rVert_{H^{\gamma}}^{2} + \lVert b\rVert_{H^{\gamma}}^{2}))(\lVert \Lambda^{\gamma}u\rVert_{L^{2}}^{2} + \lVert \Lambda^{\gamma} b\rVert_{L^{2}}^{2})
\end{eqnarray*}

where we used the Gagliardo-Nirenberg inequality, Proposition 3.1 and Young's inequality.

Thus, we can obtain this estimate for $\gamma = 1, 2, \hdots, \lfloor 2\beta + 1\rfloor + 1$ and then sum to obtain 

\begin{eqnarray}
&&\partial_{t}(\lVert u\rVert_{H^{s}}^{2} + \lVert b\rVert_{H^{s}}^{2}) + 2\eta \lVert \Lambda^{s + \beta} b\rVert_{L^{2}}^{2}\\
&\leq& c(u_{0}, b_{0}, T) (\lVert \Lambda^{2\beta}b\rVert_{L^{2}}^{2} + 1)(\log_{2}(2+ \lVert u\rVert_{H^{s}}^{2} + \lVert b\rVert_{H^{s}}^{2}))(\lVert u\rVert_{H^{s}}^{2} + \lVert b\rVert_{H^{s}}^{2})\nonumber
\end{eqnarray}

Integrating in time and using Proposition 3.2 complete the proof of Theorem 1.1.  

\section{Proof of Theorem 1.2}

We work on 

\begin{equation}
\begin{cases}
\partial_{t} u +(u\cdot\nabla) u - (b\cdot\nabla) b + \nabla \pi + \nu \Lambda^{2\alpha} u = 0\\
\partial_{t} b + (u\cdot\nabla)b - (b\cdot\nabla) u + \eta \Lambda^{2\beta} b = 0
\end{cases}
\end{equation}

Taking $L^{2}$-inner products of the first equation with $u$ and the second with $b$, using incompressibility conditions again, we obtain 

\begin{equation}
\sup_{t\in [0,T]} \left(\lVert u(t)\rVert_{L^{2}}^{2} + \lVert b(t)\rVert_{L^{2}}^{2}\right) + 2\nu\int_{0}^{T} \lVert \Lambda^{\alpha} u\rVert_{L^{2}}^{2}d\tau + 2\eta \int_{0}^{T} \lVert \Lambda^{\beta} b\rVert_{L^{2}}^{2} d\tau \leq c(u_{0}, b_{0}, T)
\end{equation}

The first proposition can be obtained similarly as before: 

\begin{proposition}
Let $N = 2, \nu, \eta > 0, \alpha \in (0, \frac{1}{2}), \beta \in (\frac{5}{4}, \frac{3}{2}]$ such that $\alpha + 2\beta > 3$. Then for any  solution pair $(u,b)$ to (1) in $[0,T]$, there exists a constant $c(u_{0}, b_{0}, T) > 0$ such that   

\begin{equation*}
\sup_{t\in [0,T]}\left(\lVert w(t)\rVert_{L^{2}}^{2} + \lVert j(t)\rVert_{L^{2}}^{2}\right) + \nu\int_{0}^{T} \lVert \Lambda^{\alpha} w\rVert_{L^{2}}^{2}d\tau + \eta \int_{0}^{T}\lVert \Lambda^{\beta} j\rVert_{L^{2}}^{2} d\tau \leq c(u_{0}, b_{0}, T)
\end{equation*}

\end{proposition}

\proof{

Taking curls on (12), we have 

\begin{eqnarray}
\partial_{t} w + \nu \Lambda^{2\alpha} w &=& - (u\cdot\nabla)w + (b\cdot\nabla)j\\
\partial_{t} j + \eta \Lambda^{2\beta}j &=& - (u\cdot\nabla)j  +  (b\cdot\nabla)w
+ 2[\partial_{1}b_{1}(\partial_{1}u_{2} + \partial_{2}u_{1}) - \partial_{1}u_{1}(\partial_{1}b_{2} + \partial_{2}b_{1})]\nonumber
\end{eqnarray}

Let us assume $\beta < \frac{3}{2}$ first. Taking $L^{2}$-inner products with $w$ and $j$ respectively, we estimate as before

\begin{eqnarray*}
\frac{1}{2} \partial_{t} X(t) + \nu \lVert \Lambda^{\alpha} w\rVert_{L^{2}}^{2} + \eta \lVert \Lambda^{\beta} j\rVert_{L^{2}}^{2} 
&\leq& c\lVert \nabla b\rVert_{L^{4}} \lVert \nabla u\rVert_{L^{2}} \lVert j\rVert_{L^{4}}\\
&\leq& c\lVert \Lambda^{\beta} b\rVert_{L^{2}}^{2\beta - 1}\lVert \Lambda^{\beta} j\rVert_{L^{2}}^{3-2\beta}\lVert w\rVert_{L^{2}}\\
&\leq& \frac{\eta}{2} \lVert \Lambda^{\beta} j\rVert_{L^{2}}^{2} + c\lVert \Lambda^{\beta} b\rVert_{L^{2}}^{2}(1+ \lVert w\rVert_{L^{2}}^{2})
\end{eqnarray*}

due to H$\ddot{o}$lder's, Gagliardo-Nirenberg and Young's inequalities.  

If $\beta = \frac{3}{2}$, by Sobolev embedding of $\dot{H}^{\frac{1}{2}}(\mathbb{R}^{2})\hookrightarrow L^{4}(\mathbb{R}^{2})$ we immediately have 

\begin{eqnarray*}
\frac{1}{2} \partial_{t} X(t) + \nu \lVert \Lambda^{\alpha} w\rVert_{L^{2}}^{2} + \eta \lVert \Lambda^{\beta} j\rVert_{L^{2}}^{2}
\leq c\lVert \nabla b\rVert_{L^{4}} \lVert \nabla u\rVert_{L^{2}} \lVert j\rVert_{L^{4}} 
\leq c\lVert \Lambda^{\beta}b\rVert_{L^{2}}^{2}(1+X(t))
\end{eqnarray*}

Absorbing the diffusive term, Gronwall's inequality completes the proof of Proposition 4.1. 

}

\begin{proposition}
Let $N = 2, \nu, \eta > 0, \alpha \in (0,\frac{1}{2}), \beta \in (\frac{5}{4}, \frac{3}{2}]$ such that $\alpha + 2\beta > 3$. Then for any solution pair $(u, b)$ to (1) and $\gamma \in (\beta, \alpha + \beta)$, there exists a constant $c(u_{0}, b_{0}, T) > 0$ such that  

\begin{equation*}
\sup_{t\in [0,T]}\lVert \Lambda^{\gamma}b(t)\rVert_{L^{2}} + \eta \int_{0}^{T} \lVert \Lambda^{\gamma + \beta} b\rVert_{L^{2}}^{2}d\tau \leq c(u_{0}, b_{0}, T)
\end{equation*}

\end{proposition}

\proof{

Let us first fix $\alpha \in (0,\frac{1}{2}), \beta \in (\frac{5}{4}, \frac{3}{2}]$ such that $\alpha + 2\beta > 3$. Then we fix $\gamma \in (\beta, \alpha + \beta)$ and estimate the second equation of (12) after multiplying by $\Lambda^{2\gamma}b$ and integrating, 

\begin{eqnarray*}
\frac{1}{2} \partial_{t} \lVert \Lambda^{\gamma} b\rVert_{L^{2}}^{2} + \eta \lVert \Lambda^{\gamma + \beta} b\rVert_{L^{2}}^{2}
&\leq& \lVert (u\cdot\nabla)b\rVert_{\dot{H}^{\gamma - \beta}}\lVert \Lambda^{\gamma + \beta} b\rVert_{L^{2}} + \lVert (b\cdot\nabla)u\rVert_{\dot{H}^{\gamma - \beta}}\lVert \Lambda^{\gamma + \beta} b\rVert_{L^{2}}\\
&\leq& \frac{\eta}{2} \lVert \Lambda^{\gamma + \beta} b\rVert_{L^{2}}^{2} + c( \lVert (u\cdot\nabla)b\rVert_{\dot{H}^{\gamma - \beta}}^{2} + \lVert (b\cdot\nabla) u\rVert_{\dot{H}^{\gamma - \beta}}^{2})
\end{eqnarray*}

by H$\ddot{o}$lder's and Young's inequalities. Now we estimate separately. By Lemma 2.5, we have 

\begin{equation*}
\lVert (u\cdot\nabla)b\rVert_{\dot{H}^{\gamma - \beta}}^{2} \leq c\lVert u\rVert_{\dot{H}^{2-\beta}}^{2} \lVert \nabla b\rVert_{\dot{H}^{\gamma - 1}}^{2}
\end{equation*}

Using Gagliardo-Nirenberg inequality, the bound on $\lVert b\rVert_{L^{2}}$ from (13) and the bound on $\lVert j\rVert_{L^{2}}$ from Proposition 4.1, we further bound by 

\begin{eqnarray*}
c\lVert u\rVert_{\dot{H}^{2-\beta}}^{2} \lVert \nabla b\rVert_{\dot{H}^{\gamma - 1}}^{2}\leq c\lVert u\rVert_{L^{2}}^{2(\beta - 1)}\lVert \nabla u\rVert_{L^{2}}^{2(2-\beta)}\lVert b\rVert_{\dot{H}^{\gamma}}^{2}
\leq c\lVert b\rVert_{\dot{H}^{\gamma}}^{2}
\end{eqnarray*}

Similarly

\begin{eqnarray*}
\lVert (b\cdot\nabla)u\rVert_{\dot{H}^{\gamma - \beta}}^{2} 
&\leq& c\lVert b\rVert_{\dot{H}^{\gamma - \beta + 1-\alpha}}^{2} \lVert \nabla u\rVert_{\dot{H}^{\alpha}}^{2}\\
&\leq& c\lVert b\rVert_{L^{2}}^{2(\alpha + \beta - \gamma)}\lVert j\rVert_{L^{2}}^{2(1-(\alpha + \beta - \gamma))}\lVert \Lambda^{\alpha} w\rVert_{L^{2}}^{2}
\leq c\lVert \Lambda^{\alpha} w\rVert_{L^{2}}^{2}
\end{eqnarray*}

by Lemma 2.5, Gagliardo-Nirenberg inequality, (13) and Proposition 4.1. 

Therefore, absorbing the diffusive term, we have shown

\begin{equation*}
\partial_{t} \lVert \Lambda^{\gamma} b\rVert_{L^{2}}^{2} + \eta \lVert \Lambda^{\gamma + \beta}b\rVert_{L^{2}}^{2}
\leq c(\lVert b\rVert_{\dot{H}^{\gamma}}^{2} + \lVert \Lambda^{\alpha} w\rVert_{L^{2}}^{2})
\end{equation*}

Integrating in time and using Proposition 4.1 allows us to complete the proof of Proposition 4.2. 

}

\begin{proposition}
Let $N = 2, \nu, \eta > 0, \alpha \in (0,\frac{1}{2}), \beta \in (\frac{5}{4}, \frac{3}{2}]$ such that $\alpha + 2\beta > 3$. Then for any solution pair $(u, b)$ to (1) in $[0,T]$, there exists a constant $c(u_{0}, b_{0}, T) > 0$ such that  

\begin{equation*}
\sup_{t\in [0,T]}\lVert w(t)\rVert_{L^{\infty}}\leq c(u_{0}, b_{0}, T)
\end{equation*}

\end{proposition}

\proof{ 

We fix $p > 2$ and also $\alpha \in (0,\frac{1}{2}), \beta \in (\frac{5}{4}, \frac{3}{2}]$ so that $\alpha + 2\beta > 3$. Then we may find $\gamma \in (\beta, \alpha + \beta)$ so that 

\begin{equation}
\gamma + \beta > 3
\end{equation}

Now we multiply the first equation of (14) by $\lvert w\rvert^{p-2}w$, integrate in space to obtain 

\begin{equation*}
\frac{1}{p} \partial_{t} \lVert w\rVert_{L^{p}}^{p} + \nu \int \Lambda^{2\alpha} w\lvert w\rvert^{p-2} w = -\int (u\cdot\nabla)w\lvert w\rvert^{p-2} w + \int (b\cdot\nabla) j \lvert w\rvert^{p-2} w
\end{equation*}

By Lemma 2.4, using incompressibility condition and H$\ddot{o}$lder's inequality we obtain

\begin{eqnarray*}
\frac{1}{p} \partial_{t} \lVert w\rVert_{L^{p}}^{p} 
\leq \int (b\cdot\nabla)j\lvert w\rvert^{p-2} w \leq \lVert b\rVert_{L^{p}} \lVert \nabla j\rVert_{L^{\infty}} \lVert w\rVert_{L^{p}}^{p-1}
\end{eqnarray*}

Dividing by $\lVert w\rVert_{L^{p}}^{p-1}$, we obtain 

\begin{eqnarray*}
\partial_{t} \lVert w\rVert_{L^{p}} 
\leq c(p) c\lVert b\rVert_{L^{2}}^{\frac{2+p(\gamma -1)}{p\gamma}}\lVert \Lambda^{\gamma} b\rVert_{L^{2}}^{\frac{p-2}{p\gamma}}\lVert \Lambda^{\gamma} b\rVert_{L^{2}}^{\frac{\gamma + \beta - 3}{\beta}}\lVert \Lambda^{\gamma + \beta}b\rVert_{L^{2}}^{\frac{3-\gamma}{\beta}}
\end{eqnarray*}

due to Gagliardo-Nirenberg inequalities that lead to  

\begin{equation*}
\lVert b\rVert_{L^{p}} \leq c(p) \lVert b\rVert_{L^{2}}^{\frac{2+p(\gamma -1)}{p\gamma}}\lVert \Lambda^{\gamma} b\rVert_{L^{2}}^{\frac{p-2}{p\gamma}}
\end{equation*}

and 

\begin{eqnarray*}
\lVert \nabla j\rVert_{L^{\infty}} \leq c\lVert \Lambda^{\gamma - 1}j\rVert_{L^{2}}^{\frac{\gamma + \beta - 3}{\beta}}\lVert \Lambda^{\gamma + \beta - 1}j\rVert_{L^{2}}^{\frac{3-\gamma}{\beta}} \leq c\lVert \Lambda^{\gamma} b\rVert_{L^{2}}^{\frac{\gamma + \beta - 3}{\beta}}\lVert \Lambda^{\gamma + \beta}b\rVert_{L^{2}}^{\frac{3-\gamma}{\beta}}
\end{eqnarray*}

where we used (15). Therefore, using the bound on $\lVert b\rVert_{L^{2}}$ from (13) and the bound on $\lVert \Lambda^{\gamma} b\rVert_{L^{2}}$ from Proposition 4.2 and Young's inequality we obtain

\begin{eqnarray*}
\partial_{t} \lVert w\rVert_{L^{p}} 
\leq  c(p) c \lVert \Lambda^{\gamma + \beta} b\rVert_{L^{2}}^{\frac{3-\gamma}{\beta}} \leq c(p) c (1+ \lVert \Lambda^{\gamma + \beta}b\rVert_{L^{2}}^{2})
\end{eqnarray*}

Integrating in time and taking limit $p\to \infty$, we have (cf. (7))

\begin{equation*}
\sup_{t \in [0,T]} \lVert w(t)\rVert_{L^{\infty}} \leq \lVert w(0)\rVert_{L^{\infty}} + c_{0}c \int_{0}^{T} 1+ \lVert \Lambda^{\gamma + \beta} b\rVert_{L^{2}}^{2} d\tau \leq c(u_{0}, b_{0}, T)
\end{equation*}

due to Proposition 4.2. This completes the proof of Proposition 4.3. 

}

\textit{Proof of Theorem 1.2}

How Proposition 4.3 leads to the higher regularity is very similar to the proof of Theorem 1.1. We sketch it for completeness. An application of Lemma 2.3 and Proposition 4.3 leads to the same bound of $\lVert \nabla u\rVert_{L^{\infty}}$ as (9). For any $\gamma > 2$, we apply $\Lambda^{\gamma}$ on (12), take $L^{2}$-inner products with $\Lambda^{\gamma}u$ and $\Lambda^{\gamma}b$ respectively to estimate using Lemma 2.2

\begin{eqnarray*}
&&\frac{1}{2}\partial_{t} (\lVert \Lambda^{\gamma} u\rVert_{L^{2}}^{2} + \lVert \Lambda^{\gamma} b\rVert_{L^{2}}^{2}) + \nu \lVert \Lambda^{\gamma+\alpha} u\rVert_{L^{2}}^{2} + \eta \lVert \Lambda^{\gamma+\beta} b\rVert_{L^{2}}^{2}\\
&\leq& c(\lVert \nabla u\rVert_{L^{\infty}} + \lVert \nabla b\rVert_{L^{\infty}}) (\lVert \Lambda^{\gamma} u\rVert_{L^{2}}^{2} + \lVert \Lambda^{\gamma} b\rVert_{L^{2}}^{2})\\
&\leq& c\left(c(u_{0}, b_{0}< T) \log_{2}(2+ \lVert u\rVert_{H^{\gamma}} + 1) + \lVert \Lambda^{\gamma} b\rVert_{L^{2}}^{\frac{\gamma + \beta - 2}{\beta}}\lVert \Lambda^{\gamma + \beta} b\rVert_{L^{2}}^{1-\frac{\gamma + \beta - 2}{\beta}}\right) (\lVert \Lambda^{\gamma} u\rVert_{L^{2}}^{2} + \lVert \Lambda^{\gamma} b\rVert_{L^{2}}^{2})\\
&\leq& c(u_{0}, b_{0}, T) (1+ \lVert \Lambda^{\gamma + \beta} b\rVert_{L^{2}}^{2}) (\log_{2}(2+ \lVert u\rVert_{H^{\gamma}}^{2} + \lVert b\rVert_{H^{\gamma}}^{2}) ) (\lVert \Lambda^{\gamma} u\rVert_{L^{2}}^{2} + \lVert \Lambda^{\gamma} b\rVert_{L^{2}}^{2})
\end{eqnarray*}

where we used (9), Gagliardo-Nirenberg and Young's inequalities. We sum over $\gamma = 1, 2, \hdots, \lfloor 2\beta + 1 \rfloor + 1$, integrate in time and use the bound of 

\begin{equation*}
\int_{0}^{T} \lVert \Lambda^{\gamma + \beta} b\rVert_{L^{2}}^{2} d\tau \leq c(u_{0}, b_{0}, T)
\end{equation*}

from Proposition 4.2 to complete the proof of Theorem 1.2. 

\section{Appendix}

\subsection{The case $\nu, \eta > 0, \alpha \geq \frac{1}{2}, \beta \geq 1$}

As remarked, the following result is immediate from our work. 

\begin{theorem}
(cf. [22])
Let $N = 2, \nu, \eta > 0, \alpha \geq \frac{1}{2}, \beta \geq 1$. Then for all initial data pair $(u_{0}, b_{0}) \in H^{s}(\mathbb{R}^{2}) \times H^{s}(\mathbb{R}^{2}), s > 2, s \geq \max\{1+2\alpha, 1+2\beta\}$, there exists a unique global strong solution pair $(u,b)$ to (1) such that 

\begin{eqnarray*}
&&u \in C([0,\infty); H^{s}(\mathbb{R}^{2}))\cap L^{2}([0,\infty); H^{s+\alpha} (\mathbb{R}^{2}))\\
&&b \in C([0,\infty); H^{s}(\mathbb{R}^{2}))\cap L^{2}([0,\infty); H^{s+\beta}(\mathbb{R}^{2}))
\end{eqnarray*}

\end{theorem}

Because our proof is simple, we sketch it here. We work on 

\begin{equation}
\begin{cases}
\partial_{t}u + (u\cdot\nabla)u - (b\cdot\nabla) b + \nabla \pi + \nu \Lambda u= 0\\
\partial_{t}b + (u\cdot\nabla) b - (b\cdot\nabla) u + \eta \Lambda^{2}b = 0
\end{cases}
\end{equation}

The following can be immediately obtained as before:

\begin{eqnarray}
\sup_{t\in [0,T]}\left(\lVert u(t) \rVert_{L^{2}}^{2} + \lVert b(t)\rVert_{L^{2}}^{2}\right) + 2\nu\int_{0}^{T} \lVert \Lambda^{\frac{1}{2}} u\rVert_{L^{2}}^{2}d\tau + 2\eta \int_{0}^{T}\lVert \Lambda b\rVert_{L^{2}}^{2} d\tau \leq c(u_{0}, b_{0})
\end{eqnarray}

Using (17), the following can be obtained as before as well: 

\begin{proposition}
Let $N = 2, \nu, \eta > 0, \alpha = \frac{1}{2}, \beta = 1$. Then for any solution pair $(u,b)$ to (1) in $[0,T]$ there exists a constant $c(u_{0}, b_{0}, T)$ such that   

\begin{equation*}
\sup_{t\in [0,T]}\left(\lVert w(t)\rVert_{L^{2}}^{2} + \lVert j(t)\rVert_{L^{2}}^{2}\right) + \nu\int_{0}^{T}\lVert \Lambda^{\frac{1}{2}} w\rVert_{L^{2}}^{2}d\tau + \eta\int_{0}^{T}\lVert \Lambda j\rVert_{L^{2}}^{2} d\tau \leq c(u_{0}, b_{0}, T)
\end{equation*}

\end{proposition}

Obtaining higher estimate from Propositions 5.2 is immediate. Applying  $\Lambda^{3}$ on the first and second equations of (16), taking $L^{2}$-inner products with $\Lambda^{3}u$ and $\Lambda^{3}b$ respectively using Lemma 2.2 one can immediately obtain 

\begin{eqnarray}
&&\frac{1}{2}\partial_{t} (\lVert \Lambda^{3} u\rVert_{L^{2}}^{2} + \lVert \Lambda^{3}b\rVert_{L^{2}}^{2}) + \nu \lVert \Lambda^{3+\frac{1}{2}}u\rVert_{L^{2}}^{2} + \eta \lVert \Lambda^{4} b\rVert_{L^{2}}^{2}\\
&\leq& \frac{\nu}{2} \lVert \Lambda^{3+\frac{1}{2}} u\rVert_{L^{2}}^{2} + \frac{\eta}{2} \lVert \Lambda^{4} b\rVert_{L^{2}}^{2} + c(\lVert \Lambda^{1+\frac{1}{2}} u\rVert_{L^{2}}^{2} + \lVert \Lambda^{2} b\rVert_{L^{2}}^{2} + 1) (\lVert \Lambda^{3} u\rVert_{L^{2}}^{2} + \lVert \Lambda^{3} b\rVert_{L^{2}}^{2})\nonumber
\end{eqnarray}

in particular by using Sobolev embedding of $\dot{H}^{\frac{1}{2}}(\mathbb{R}^{2}) \hookrightarrow L^{4}(\mathbb{R}^{2})$. This leads to the completion of the proof of Theorem 5.1.

\subsection{Proof of Lemma 2.5}

In this subsection, for readers' convenience, we sketch the proof of Lemma 2.5. Let us recall the notion of Besov spaces (cf. [5]). We denote by $\mathcal{S}(\mathbb{R}^{2})$ the Schwartz class functions and $\mathcal{S}'(\mathbb{R}^{2})$, its dual. We define $\mathcal{S}_{0}$ to be the subspace of $\mathcal{S}$ in the following sense:

\begin{equation*}
\mathcal{S}_{0} = \{\phi \in \mathcal{S}, \int_{\mathbb{R}^{2}}\phi(x)x^{\gamma}dx = 0, \lvert\gamma\rvert = 0, 1, 2, ... \}
\end{equation*}

Its dual $\mathcal{S}_{0}'$ is given by $\mathcal{S}_{0}' = \mathcal{S}/\mathcal{S}_{0}^{\perp} = \mathcal{S}'/\mathcal{P}$ where $\mathcal{P}$ is the space of polynomials. For $j \in \mathbb{Z}$ we define 

\begin{equation*}
A_{j} = \{\xi\in \mathbb{R}^{2}: 2^{j-1} < \lvert\xi\rvert < 2^{j+1}\}
\end{equation*}

It is well-known that there exists a sequence $\{\Phi_{j}\} \in \mathcal{S}(\mathbb{R}^{2})$ such that 

\begin{equation*}
\text{ supp }\hat{\Phi}_{j}\subset A_{j}, \hspace{5mm} \hat{\Phi}_{j}(\xi) = \hat{\Phi}_{0}(2^{-j}\xi) \hspace{5mm} \text{or} \hspace{5mm} \Phi_{j}(x) = 2^{2j}\Phi_{0}(2^{j}x) \hspace{5mm} \text{ and } 
\end{equation*}

\begin{equation*}
\sum_{j = -\infty}^{\infty}\hat{\Phi}_{j}(\xi) = 
\begin{cases}
1 \hspace{2mm} \text{ if } \hspace{1mm} \xi \in \mathbb{R}^{2}\setminus \{0\}\\
0 \hspace{2mm} \text{ if } \hspace{1mm} \xi = 0
\end{cases}
\end{equation*}

To define the homogeneous Besov space, we set 

\begin{equation*}
\dot{\Delta}_{j}f = \Phi_{j}\ast f, \hspace{5mm} j = 0, \pm 1, \pm 2, ...
\end{equation*}

With such we can define for $s \in \mathbb{R}, p, q \in [1,\infty]$, the homogeneous Besov space 

\begin{equation*}
\dot{B}_{p, q}^{s} = \{f \in \mathcal{S}_{0}' : \lVert f\rVert_{\dot{B}_{p, q}^{s}} < \infty\}
\end{equation*}

where 

\begin{equation*}
\lVert f\rVert_{\dot{B}_{p,q}^{s}} = 
\begin{cases}
\left(\sum_{j}(2^{js}\lVert\dot{\Delta}_{j}f\rVert_{L^{p}})^{q}\right)^{\frac{1}{q}} & \text{if } q < \infty\\
\sup_{j}2^{js}\lVert\dot{\Delta}_{j}f\rVert_{L^{p}} & \text{if } q = \infty
\end{cases}
\end{equation*}

To define the inhomogeneous Besov space, we let $\Psi \in C_{0}^{\infty}(\mathbb{R}^{2})$ be such that

\begin{equation*}
1 = \hat{\Psi}(\xi) + \sum_{j=0}^{\infty}\hat{\Phi}_{j}(\xi), \hspace{5mm} \Psi\ast f + \sum_{j=0}^{\infty}\Phi_{j}\ast f = f
\end{equation*}

for any $f \in \mathcal{S}'$. With that, we set

\begin{equation*}
\Delta_{j}f = 
\begin{cases}
0 &\text{ if } j \leq -2\\
\Psi\ast f, &\text{ if } j = -1\\
\Phi_{j}\ast f, &\text{ if } j = 0, 1, 2, ...
\end{cases}
\end{equation*}

and define for any $s \in \mathbb{R}, p, q \in [1, \infty]$, the  inhomogeneous Besov space 

\begin{equation*}
B_{p, q}^{s} = \{f \in \mathcal{S}' : \lVert f\rVert_{B_{p, q}^{s}} < \infty \}
\end{equation*}

where

\begin{equation*}
\lVert f\rVert_{B_{p, q}^{s}} = 
\begin{cases}
\left(\sum_{j = -1}^{\infty}(2^{js}\lVert\Delta_{j}f\rVert_{L^{p}})^{q}\right)^{\frac{1}{q}}, &\text{if } q < \infty\\
\sup_{-1 \leq j < \infty}2^{js}\lVert \Delta_{j}f\rVert_{L^{p}} &\text{if } q = \infty
\end{cases}
\end{equation*}

In particular $\dot{B}_{2,2}^{s} = \dot{H}^{s}, B_{2,2}^{s} = H^{s}$. The following lemma is very useful:

\begin{lemma} (cf. [5])
Bernstein's Inequality: Let f $\in L^{p}(\mathbb{R}^{2})$ with 1 $\leq p \leq q \leq \infty$ and $0 < r < R$. Then for all $k \in \mathbb{Z}^{+}\cup\{0\}$, and $\lambda > 0$, there exists a constant $C_{k} > 0$ such that

\begin{equation*}
\begin{cases}
\sup_{\lvert \gamma \rvert = k} \lVert\partial^{\gamma}f\rVert_{L^{q}} \leq C_{k}\lambda^{k+2(\frac{1}{p} - \frac{1}{q})}\lVert f\rVert_{L^{p}}  & \text{ if } \text{supp }\hat{f} \subset \{\xi : \lvert\xi\rvert \leq \lambda r\}\\
C_{k}^{-1}\lambda^{k}\lVert f \rVert_{L^{p}} \leq sup_{\lvert\gamma\rvert = k} \lVert\partial^{\gamma}f\rVert_{L^{p}} \leq C_{k}\lambda^{k}\lVert f \rVert_{L^{p}}  & \text{ if  supp } \hat{f} \subset \{\xi : \lambda r \leq \lvert\xi\rvert \leq \lambda R\}
\end{cases}
\end{equation*}

and if we replace derivative $\partial^{\gamma}$ by the fractional derivative, the inequalities remain valid only with trivial modifications.
\end{lemma}

Now recall Bony's paraproduct decomposition (cf. [5]): 

\begin{equation*}
fg = \dot{T}(f,g) + \dot{R}(f,g) + \dot{T}(g,f)
\end{equation*}

where with $\dot{S}_{j}f = \sum_{l\leq j-1}\dot{\Delta}_{l}f$, 

\begin{equation*}
\dot{T}(f,g) = \sum_{j} \dot{S}_{j-1} f \dot{\Delta}_{j}g, \hspace{5mm} \dot{R}(f,g) = \sum_{i=-1}^{1} \sum_{j} \dot{\Delta}_{j} f \dot{\Delta}_{j+i} g, \hspace{5mm} \dot{T}(g,f) = \sum_{j}\dot{S}_{j-1} g \dot{\Delta}_{j} f
\end{equation*}

On the estimate of $\dot{T}(f,g)$, we make use of the hypothesis that  $\sigma_{1} < 1$, on the estimate of $\dot{T}(g,f)$ that $\sigma_{2} < 1$, and on the estimate of $\dot{R}(f,g)$ that $\sigma_{1} + \sigma_{2} > 0$. Firstly,

\begin{eqnarray*}
\lVert \dot{\Delta}_{k} \dot{T}(f,g) \rVert_{L^{2}}
&\leq& c \sum_{j: \lvert j-k\rvert \leq 1}\sum_{l \leq j-2}\lVert \dot{\Delta}_{l} f\dot{\Delta}_{j} g\rVert_{L^{2}}\\ 
&\leq& c \sum_{j: \lvert j-k\rvert \leq 1}\sum_{l\leq j-2}\lVert \dot{\Delta}_{l} f\rVert_{L^{\infty}} \lVert \dot{\Delta}_{j} g\rVert_{L^{2}}\\
&\leq& c\sum_{j: \lvert j-k\rvert \leq 1}\sum_{l \leq j-2}2^{(l-j)(1-\sigma_{1})}2^{l\sigma_{1} + j - j\sigma_{1}}\lVert \dot{\Delta}_{l} f\rVert_{L^{2}} \lVert \dot{\Delta}_{j} g\rVert_{L^{2}}\\
&\leq& c\sum_{j: \lvert j-k\rvert \leq 1}\lVert f\rVert_{\dot{B}_{2, \infty}^{\sigma_{1}}}2^{j(1-\sigma_{1})}\lVert \dot{\Delta}_{j} g\rVert_{L^{2}}
\end{eqnarray*}

by H$\ddot{o}$lder's and Bernstein's inequalities. Therefore, multiplying by $2^{k(\sigma_{1} + \sigma_{2} - 1)}$ and taking $l^{2}$-norm, we have 

\begin{eqnarray*}
\lVert 2^{k(\sigma_{1} + \sigma_{2} - 1)}\lVert \dot{\Delta}_{k} \dot{T}(f,g)\rVert_{L^{2}} \rVert_{l^{2}}&\leq& c\lVert f\rVert_{\dot{B}_{2,\infty}^{\sigma_{1}}}\lVert \sum_{j: \lvert j-k\rvert \leq 1}2^{(j-k)(1-\sigma_{1})}2^{k\sigma_{2} }\lVert \dot{\Delta}_{j} g\rVert_{L^{2}} \rVert_{l^{2}}\\
&\leq& c\lVert f\rVert_{\dot{B}_{2,\infty}^{\sigma_{1}}}\lVert 2^{k\sigma_{2}}\lVert \dot{\Delta}_{k} g\rVert_{L^{2}} \rVert_{l^{2}} \leq c\lVert f\rVert_{\dot{B}_{2,\infty}^{\sigma_{1}}}\lVert g\rVert_{\dot{B}_{2,2}^{\sigma_{2}}}
\end{eqnarray*}

where we used that because $\lvert j-k\rvert \leq 1$, we may replace $j$ by $k$ modifying constant. Similarly, using H$\ddot{o}$lder's and Bernstein's inequalities, one can show 

\begin{eqnarray*}
\lVert \dot{\Delta}_{k} \dot{T}(g,f)\rVert_{L^{2}}
&\leq& c\sum_{j: \lvert j-k\rvert \leq 1} \sum_{l \leq j-2} \lVert \dot{\Delta}_{l} g\rVert_{L^{\infty}} \lVert \dot{\Delta}_{j} f\rVert_{L^{2}}\\
&\leq& c\sum_{j: \lvert j-k\rvert \leq 1} \sum_{l \leq j-2} 2^{(l-j)(1-\sigma_{2})}2^{l\sigma_{2} + j(1-\sigma_{2})}\lVert \dot{\Delta}_{l} g\rVert_{L^{2}} \lVert \dot{\Delta}_{j} f\rVert_{L^{2}}\\
&\leq& c\sum_{j: \lvert j-k\rvert \leq 1}\lVert g\rVert_{\dot{B}_{2,\infty}^{\sigma_{2}}}2^{j(1-\sigma_{2})}\lVert \dot{\Delta}_{j} f\rVert_{L^{2}}
\end{eqnarray*}

This gives 

\begin{eqnarray*}
\lVert 2^{k(\sigma_{1} + \sigma_{2} - 1)}\lVert \dot{\Delta}_{k} \dot{T}(g,f)\rVert_{L^{2}} \rVert_{l^{2}}
&\leq& c\lVert g\rVert_{\dot{B}_{2,\infty}^{\sigma_{2}}}\lVert \sum_{j: \lvert j-k\rvert \leq 1}2^{(j-k)(1-\sigma_{2})}2^{k\sigma_{1}}\lVert \dot{\Delta}_{j} f\rVert_{L^{2}}\rVert_{l^{2}}\\
&\leq&  c \lVert g\rVert_{\dot{B}_{2,\infty}^{\sigma_{2}}}\lVert f\rVert_{\dot{B}_{2,2}^{\sigma_{1}}}
\end{eqnarray*}

Finally, 

\begin{eqnarray*}
\lVert \dot{\Delta}_{k} \dot{R}(f,g)\rVert_{L^{2}}
&\leq& c2^{k} \sum_{i=-1}^{1} \sum_{j: k < j} \lVert \dot{\Delta}_{j} f \dot{\Delta}_{j+i} g\rVert_{L^{1}}\\
&\leq& c \sum_{i=-1}^{1} \sum_{j: k < j} 2^{(k-j)(\sigma_{1} + \sigma_{2})}2^{k+(j-k)(\sigma_{1} + \sigma_{2})}\lVert \dot{\Delta}_{j} f\rVert_{L^{2}} \lVert \dot{\Delta}_{j+i} g\rVert_{L^{2}}
\end{eqnarray*}

by Bernstein's and H$\ddot{o}$lder's inequalities. Therefore, 

\begin{eqnarray*}
&& \lVert 2^{k(\sigma_{1} + \sigma_{2} - 1)}\lVert \dot{\Delta}_{k} \dot{R}(f,g) \rVert_{L^{2}} \rVert_{l^{2}}\\
&\leq& c\sum_{i=-1}^{1}\lVert 2^{-\lvert k\rvert (\sigma_{1} + \sigma_{2})}\rVert_{l^{2}} \lVert 2^{k(\sigma_{1} + \sigma_{2})}\lVert \dot{\Delta}_{k} f\rVert_{L^{2}} \lVert \dot{\Delta}_{k+i} g\rVert_{L^{2}} \rVert_{l^{1}}\\
&\leq& c\sum_{i=-1}^{1} \lVert 2^{k\sigma_{1}}\lVert \dot{\Delta}_{k} f\rVert_{L^{2}} \rVert_{l^{2}} \lVert 2^{k\sigma_{2}}\lVert \dot{\Delta}_{k+i} g\rVert_{L^{2}} \rVert_{l^{2}}\\
&\leq& c \lVert f\rVert_{\dot{H}^{\sigma_{1}}} \lVert g\rVert_{\dot{H}^{\sigma_{2}}}
\end{eqnarray*}

by Young's inequality for convolution and H$\ddot{o}$lder's inequality. 

\subsection{Proof of Lemma 2.3}

We fix $s > 2$ and estimate 

\begin{eqnarray*}
\lVert \nabla f\rVert_{L^{\infty}}
&\leq& c(\lVert f\rVert_{L^{2}} + \sum_{j=0}^{n-1} \lVert \Delta_{j} \nabla f\rVert_{L^{\infty}} + \sum_{j=n}^{\infty} \lVert \Delta_{j} \nabla f\rVert_{L^{\infty}})\\
&\leq& c (\lVert f\rVert_{L^{2}} + \lVert \text{curl } f\rVert_{B_{\infty, \infty}^{0}}n + \sum_{j=n}^{\infty} 2^{j(2-s)}2^{j(s-1)}\lVert \Delta_{j} \nabla f\rVert_{L^{2}})\\
&\leq& c (\lVert f\rVert_{L^{2}} + \lVert \text{curl } f\rVert_{L^{\infty}} n + 2^{n(2-s)}\lVert f\rVert_{H^{s}})
\end{eqnarray*}

by Young's inequality for convolution, Bernstein's inequalities and continuity of Riesz transform. Choosing  $n = \left(\frac{1}{s-2}\right) \log_{2}(2+ \lVert u\rVert_{H^{s}})$ immediately implies the desired result.

\end{document}